\documentclass[12pt]{article}
\usepackage{amsmath,amssymb,times}

\usepackage[dvips]{graphicx}
\usepackage{color}
   \newcommand{\blue}[1]{\textcolor{blue}{#1}}

\usepackage{dcolumn}
\usepackage{bm}

\begin{document}
\pagestyle{plain}

\title{One-Dimensional Birth-Death Process and 
Delbr\"{u}ck-Gillespie Theory of Mesoscopic
Nonlinear Chemical Reactions}

\author{Yunxin Zhang$^1$\footnote{xyz@fudan.edu.cn},\ 
Hao Ge$^{2,1}$\footnote{haoge@pku.edu.cn}\ \ and\ 
Hong Qian$^{3,1,4}$\footnote{hqian@u.washington.edu}\\[20pt]
$^1$School of Mathematical Sciences\\
and Centre for Computational Systems Biology\\ 
Fudan University, Shanghai 200433 PRC\\[10pt]
$^2$Beijing International Center for Mathematical 
Research\\ 
and Biodynamic Optical Imaging Center\\
Peking University, Beijing 100871 PRC\\[10pt]
$^3$Department of Applied Mathematics\\ 
University of Washington, Seattle, WA 98195, USA\\[10pt]
$^4$College of Mathematics, Jilin University\\
Changchun, 130012 PRC. 
}

\maketitle

\newpage

\begin{abstract}
As a mathematical theory for the stochasstic, nonlinear
dynamics of individuals within a population,  
Delbr\"{u}ck-Gillespie process (DGP) $n(t)\in\mathbb{Z}^N$,
is a birth-death system with state-dependent rates which 
contain the system size $V$ as a natural parameter.  For 
large $V$, it is intimately related to an autonomous, 
nonlinear ordinary differential equation as well as a 
diffusion process.  For nonlinear dynamical systems with 
multiple attractors, the quasi-stationary and stationary
behavior of such a birth-death process can be 
underestood in terms of a separation of time scales by a
$T^*\sim e^{\alpha V}$ $(\alpha>0)$: a relatively 
fast, intra-basin diffusion for
$t\ll T^*$ and a much slower inter-basin Markov 
jump process for $t\gg T^*$.  In the present paper for
one-dimensional systems, we study both stationary 
behavior ($t=\infty$) in terms of invariant 
distribution $p_n^{ss}(V)$, and finite time 
dynamics in terms of the mean first passsage time 
(MFPT) $T_{n_1\rightarrow n_2}(V)$.  We obtain an 
asymptotic expression of MFPT in terms of the ``stochastic
potential'' $\Phi(x,V)=-(1/V)\ln p^{ss}_{xV}(V)$.
We show in general no continuous diffusion process 
can provide asymptotically accurate representations 
for both the MFPT and the $p_n^{ss}(V)$ for a DGP.  
When $n_1$ and $n_2$ belong to two different basins of 
attraction, the MFPT yields the $T^*(V)$ 
in terms of $\Phi(x,V)\approx \phi_0(x)+(1/V)\phi_1(x)$.
For systems with a saddle-node bifurcation and catastrophe, discontinuous ``phase transition'' emerges, which can
be characterized by $\Phi(x,V)$ in the limit of $V\rightarrow\infty$.   In terms of time scale separation, the relation between deterministic, local nonlinear bifurcations
and stochastic global phase transition is discussed.  The 
one-dimensional theory is a pedagogic first step toward a 
general theory of DGP.
\end{abstract}




\section{Introduction}

Nonlinear ordinary differential equations (ODEs) and 
diffusion processes are two important mathematical 
models, respectively, for dynamics of deterministic and 
stochastic systems.  To understand the mathematical 
properties of these dynamical models, it is obligatory to 
first have a thorough analysis of one-dimensional (1-d) 
systems.  In the case of a nonlinear ODE, this is 
\begin{subequations}
\begin{equation}
		\frac{dx(t)}{dt}=b(x),
\label{the_ode}
\end{equation} 
where $x(t)$ is the state of a system at time $t$,
and in the case of diffusion processes, it is
\begin{equation}
	\frac{\partial u(x,t)}{\partial t}
	= \frac{\partial}{\partial x}\left(
	D(x)\frac{\partial u(x,t)}{\partial x}
	- b(x) u(x,t)\right),
\label{1dd}
\end{equation}
in which $u(x,t)$ is the probability density for
a system being in state $x$ at time $t$.
A wealth of mathematics has been created by
thorough investigations of these simple systems.
They are now textbook materials with great
pedagogic values \cite{perko,omalley,mazo,schuss}.
When $D(x)=0$, Eq. (\ref{1dd}) is reduced to 
Eq. (\ref{the_ode}) via the method of characteristics;
(\ref{1dd}) is known as the Liouville equation 
of (\ref{the_ode}) in phase space.
The solution to Eq. (\ref{1dd}) with vanishing
$D(x)$ can be considered as a {\em viscosity solution}
to the first-order, hyperbolic partial differential
equation.

	In recent years, in connection to mesoscopic size,
cellular biochemical dynamics, a new
type of mathematical models has emerged: the 
multi-dimensional birth-death process. An 
$N$-dimensional birth-death process is a Markov 
jump process with discrete state $\vec{\ell}\in \mathbb{Z}^N$ 
and continuous time $t$ \cite{taylor_karlin}.  When applied 
to nonlinear biochemical reaction systems \cite{epstein}, 
its time-dependent probability mass distribution, 
$p_{\vec{\ell}}(t)$ 
satisfies the {\em Chemical Master Equation} (CME), first
studied by M. Delbr\"{u}ck, while its stochastic 
trajectories can be sampled according to the Gillespie 
algorithm \cite{dgp_1,dgp_2,qian_review_1,
qian_review_2,qian_review_3}.

	The new theory for the Markov dynamics of population 
systems deserves more attentions from applied
mathematicians \cite{kurtz_book}.  In addition to
its own importance in applications, it also provides a 
unique opportunity for studying the relationship between
dynamics at mesoscopic and macroscopic levels, which in
the past has been studied mainly in terms of diffusion 
processes with Brownian noise.  It is a widely 
hold belief that birth-death processes
can be approximated by diffusions.  This turns out not
to be the case for nonlinear systems with multiple 
attractors, as we shall show.

	With this backdrop in mind, it is again obligatory 
to first carry out an comprehensive analysis for a 
1-d CME system.  Doering {\em et al.} have 
conducted an extensive investigation for the asymptotic 
expressions of the mean first passage time (MFPT) \cite{doering_1,doering_2}.  The aim of the present 
work is not on this {\em per se}, but to illustrate the 
overall mathematical structure of stochastic nonlinear 
population dynamics in terms of the 1-d system.

	The  CME for a 1-d birth-death 
process takes the form
\begin{equation}
	\frac{d}{dt}p_n(t) = u_{n-1}p_{n-1}(t)
		-(w_n+u_n)p_n(t) +  w_{n+1}p_{n+1}(t),	
		\ \ (n\ge 0)
\label{1dbdp}
\end{equation}
\end{subequations}
in which state-dependent birth and death rates
$u_n(V)$ and $w_n(V)$ are in general functions of $n$
as well as a crucial parameter $V$, the spatial size
or any other extensive quantity of 
the reaction system.  For chemical systems consisting 
of only first-order, linear reactions, both $u_n$ and 
$w_n$ are independent of $V$ \cite{othmer_05,qian_jcp_06,
huisinga_07}.  Linear systems have found wide 
applications in modeling stochastic dynamics of 
single biological macromolecules \cite{qian_review_1}, 
such as in single-molecule enzymology and molecular-motor
chemomechanics \cite{qian_bpc_02,qian_bj_08, zhang_fisher_1,zhang_fisher_2}.

	The dependence on $V$ gives rise to a very special
feature of the theory of Delbr\"{u}ck-Gillespie processes 
(DGP) and its corresponding CME: One can study the important
relation between a stochastic dynamical model with a small 
$V$ and a nonlinear deterministic dynamical system with 
infinitely large $V$ \cite{kurtz,keizer}. 
T.G. Kurtz's theorem precisely establishes such a convergence
from the stochastic trajectories of a DGP
 to the solution of a nonlinear 
ODE like Eq. (\ref{the_ode}).  In the 1-d DGP, each of 
(\ref{the_ode}), (\ref{1dd}) and (\ref{1dbdp}) 
has a role.  There is also a substantial 
difference between the stochastic system (\ref{1dd}) in which
the stochasticity $D(x)$ and deterministic $b(x)$ are not
related {\em per se}, while the stochasticity is intrinsic
in the dynamics of (\ref{1dbdp}).  Therefore models based on 
diffusion processes are often {\em phenomenological}, while
the discrete model provides a more faithful representation of 
a system's emerging dynamics based on individual's stochastic 
behavior.  

\blue{Motivated by biochemical applications, recent studies
on 1-d DGPs have focused on ($i$) nonlinear 
bistability, stochastic bistability and multi-stability 
(i.e., multi-modal distribution) \cite{vellela_qian_07,vellela_qian_09,bishop_qian_bj_10,
ge_qian_chaos}, 
($ii$) non-equilibrium thermodynamics and phase transition
\cite{vellela_qian_09,ge_qian_prl,ge_qian_jrsi}, 
($iii$) large $V$ asymptotics in terms of large deviation theory
\cite{ge_qian_jrsi,ge_qian_chaos},
($iv$) van't Hoff-Arrhenius analysis motivated by classical thermodynamics \cite{zhang_ge_qian}, and 
($v$) dynamics on the circle $\mathbb{S}^1$ and oscillations in terms of a rotational random walk \cite{qian_pnas_02,vellela_qian_10,ge_qian_chaos}.
}

\blue{A series of mathematical issues arise in 
the investigations.  The present paper
initiates a systematic treatment of some of 
them.}
One might be surprised by that there
are still significant unresolved mathematical 
questions for a one-dimensional birth-death
process. We simply point out that for large $V$, 
the problem under investigation is intimately
related to the Eq. (\ref{1dd}) with a singularly 
perturbed coefficient $D(x)\propto (1/V)$.
This is still an active area of
research on its own \cite{omalley_sirev}. 
In addition, even though straightforward,  many 
explicit formulae in connection to the one-dimensional 
Eq. (\ref{1dbdp}), also known as hopping models in 
statistical physics, had not been obtained until a 
need arose from applications.  A case in point was the 1983
paper of B. Derrida \cite{derrida}. See also  
\cite{doering_1,doering_2} for recent work on 
the asymptotic analysis of the MFPT problem.  Finally, the
newly introduced van't Hoff-Arrhenius 
analysis \cite{zhang_ge_qian} and the analysis for limit
cycles \cite{ge_qian_chaos} both require consistent asymptotic 
expansions for large $V$ beyond the usual leading order.  

\begin{figure}[t]
\centerline{\includegraphics[width=3.5in]{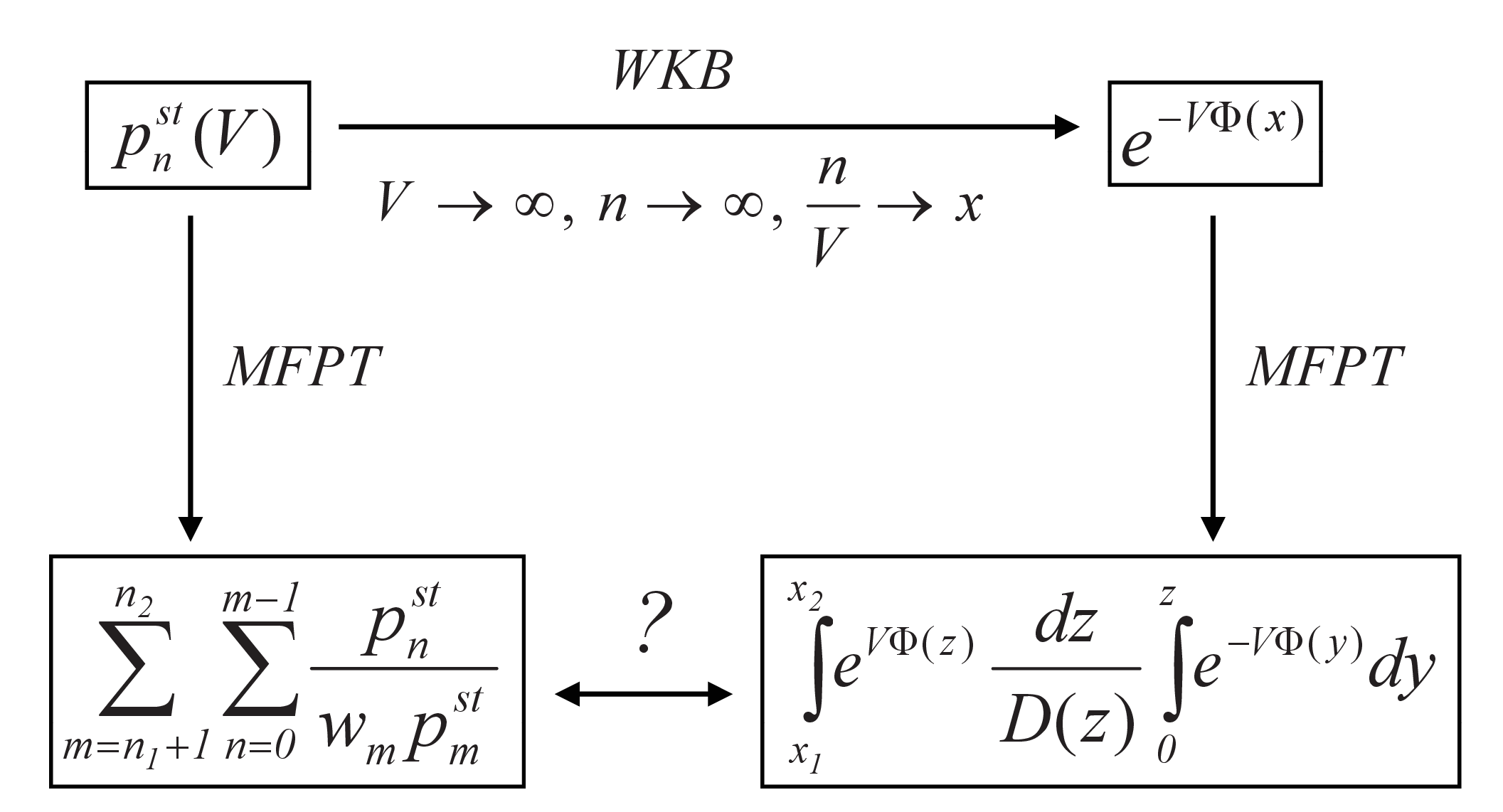}}
\caption{Logical schematics showing, for a 1-d DGP, the 
mathematical relations between infinite-time stationary 
distribution $p_n^{st}$, MFPT for finite-time dynamics, 
and their $V\rightarrow\infty$ asymptotics. MFPT for 1-d
DGP can be exactly expressed in terms of $p^{st}_n$ as
given in the lower-left box (Eq. \ref{T12_exact}); 
$p^{st}_n$ also has an asymptotic form shown in 
the upper-right box.  For 1-d continuous diffusion, its 
stationary density function is related to MFPT as shown 
by the two boxes on the right (Eq. \ref{mfpt}).  The
two MFPTs are ``analogous'' if we identify $w_m$ 
with $D(z)$ and replace summations with integrals. 
A remaining question: What is the asymptotic expression
for the MFPT in terms of the asymptotic stochasstic
potential $\Phi(x)$.
}
\label{fig_1}
\end{figure}

	One of the questions we
study in the present work can be succinctly
described in terms of the diagram in Fig. \ref{fig_1}.
It is well established that in the limit
of large $V$, the stationary solution to the
1-d (\ref{1dbdp}) has a WKB 
(Wentzel-Kramers-Brillouin) type asymptotic expansion 
$p_n(V)\sim\exp\big(-V\phi_0(x)-\phi_1(x)\big)$ where $x=n/V$
\cite{nicolis_77,hugang,schuss}.  In chemical terms,
$n$ is the copy number of a chemical species and $x$ is 
its concentration.  Furthermore, it is straightforward 
to compute the MFPTs for both discrete birth-death 
processes and continuous diffusion \cite{vanKampen,gardiner}.  
But it is  unclear, as indicated by the question mark in Fig.
\ref{fig_1}, whether and how the MFPT of the 
birth-death processes in the limit of large $V$ is 
related to the ``stochastic potential function'' $\phi_0(x)+(1/V)\phi_1(x)$ 
obtained from the WKB expansion \cite{nicolis_77,hugang}.
This is answered in Eqs. (\ref{asym_4_T}) and (\ref{new_k_eq}). 

For the stationary solution of Kolmogorov forward
equation (\ref{1dbdp}), we now have a good 
understanding: For nonlinear dynamical systems with 
two attractors, there is an exponentially large time, 
$e^{\alpha V}$ ($\alpha>0)$ that separates the 
{\em intra-basin} dynamics in terms of
Gaussian processes \cite{qian_review_3} from
{\em the inter-basin} dynamics of a 
Markov jump process between two discrete states.
The ``boundary layer'' in the singularly
perturbed problem is precisely where different 
attractors join \cite{omalley,omalley_sirev}.
$-e^{-\alpha V}$ is in fact the second largest
eigenvalue of the linear system (\ref{1dbdp}), with 
zero being the largest one.  When $V=\infty$, there 
is a breakdown of ergodicity \cite{hanggi_84,vellela_qian_07,vellela_qian_09}.

	Recognizing this exponentially 
large time is the key to resolve the so-called {\em Keizer's paradox} \cite{vellela_qian_07,keener_1,keener_2}
which illustrates the two completely different 
pictures for the ``steady states'' of a deterministic
system and its CME counterpart.  It is also the
key to understand the difficulty of approximating
a CME like (\ref{1dbdp}) with a diffusion equation like 
(\ref{1dd}) for systems with multistability \cite{hanggi_84,
vellela_qian_09,qian_review_3}.  See more discussions
below.

	The MFPT is the solution to the time-independent
backward equation with an inhomogeneous term $-1$ \cite{vanKampen,
gardiner}.  An ambiguity arises in the asymptotics of 
MPFT as a WKB solution to the backward 
equation \cite{doering_2}.  This is reminiscent of the WKB 
approach to the stationary forward equation in terms of 
the nonlinear Hamilton-Jacobi equation \cite{ge_qian_HJE}.
One of the results in the present work, however, is
the asymptotic MFPT {\em in relation} to the 
asymptotic stationary solution to the corresponding 
forward equation.

\section{Background on Diffusion Processes}

	Because of the intimate relationship between 
Eqs. (\ref{1dbdp}) and (\ref{1dd}), we shall give a brief
summary of the relevant results for one-dimensional
continuous diffusion in Sec. \ref{con_dif}.  Even though 
it contains no new mathematical result, the presentation is
novel.  Then in Sec. \ref{k_para}, we discuss
Keizer's paradox from a novel perspective by
considering a second-order correction to the Kramers-Moyal 
expansion \cite{vanKampen,gardiner}.

\subsection{Continuous diffusion}
\label{con_dif}

For a continuous diffusion process
with $\epsilon$-small diffusion coefficient:
\begin{equation}
	\frac{\partial f(x,t)}{\partial t}
	=\frac{\partial}{\partial x}
	\left(\epsilon D(x)\frac{\partial f}{\partial x}
		- b(x) f\right),  \ \ (D(x)>0)
\end{equation}
the stochastic potential function
\begin{equation}
	\Psi(x) = -\int_0^x \frac{b(z)}{D(z)}dz
\label{the_Psi}
\end{equation}
plays a central role in its dynamics.  In terms
of the $\Psi(x)$, one has the stationary
distribution
\begin{equation}
	f^{ss}(x) = A e^{-\frac{1}{\epsilon}\Psi(x)},
\end{equation}
where $A$ is a normalization factor. $\Psi(x)$ is also
a Lyapunov function of the ordinary differential
equation dynamics $\frac{dx}{dt}=b(x)$ since 
\begin{equation}
	\frac{d}{dt}\Psi(x(t)) = \frac{d\Psi(x)}{dx}b(x)
		= -\frac{b^2(x)}{D(x)} \le 0.
\end{equation}
Furthermore, the mean first passage time arriving at 
$x_2$ starting at $x_1$ with a reflecting boundary at
$x_0$ $(x_0<x_1<x_2)$ is \cite{vanKampen,gardiner}
\begin{equation}
	T_{x_1\rightarrow x_2} = \int_{x_1}^{x_2}
		e^{\Psi(z)/\epsilon}\frac{dz}{D(z)}
		\int_{x_0}^ze^{-\Psi(y)/\epsilon}dy.
\label{mfpt}
\end{equation}
On the other hand, solving the stationary flux $J$
passing through $x_2$ with Dirichlet boundary  value
$f(x_2)=0$, leaving the boundary value at $x_1$ unspecified,
but enforcing a normalization condition 
$\int_{x_0}^{x_2}f^{ss}(x)dx=1$,
\begin{equation}
	J^{-1} = \int_{x_1}^{x_2} e^{-\Psi(y)/\epsilon}dy
		\int_y^{x_2} e^{\Psi(z)/\epsilon}
				\frac{dz}{D(z)}.
\label{ssflux}
\end{equation}

	Note that Eqs. (\ref{mfpt}) and (\ref{ssflux}) are
exactly the same if $x_0=x_1$ in (\ref{mfpt}).  
To understand the origin of this intriguing
observation, consider the following {\em Gedankenexperiment}: 
Let a diffusing particle start at $x_1=x_0$, which is also a 
reflecting boundary. The particle can only move rightward,
and as soon as it hits $x_2 (>x_1)$, one 
immediately takes it back to $x_1$.  Repeating
this procedure forms a {\em renewal process}.  Then the
mean renewal time is $T_{x_1\rightarrow x_2}$ in 
Eq. (\ref{mfpt}).   Now imagine that one connects $x_2$
with $x_1$ to form a circle, and installs a one-way
permeable membrane at the $x_2$-$x_1$ junction: a particle
that hits from the $x_2$ side goes through the membrane and starts 
at $x_1$ instantaneously; but a particle that hits from the $x_1$ side
is reflected.  The stationary distribution for the 
diffusion particle then
satisfies $f^{ss}(x_2)=0$, $\int_{x_1}^{x_2} f^{ss}(x)dx=1$, 
and a constant flux $J(x_1) = J(x_2)$ is the $J$ in 
Eq. (\ref{ssflux}).  

	According to the {\em elementary renewal theorem} 
\cite{taylor_karlin}, $T_{x_1\rightarrow x_2}=J^{-1}$. 

Another problem which is widely employed in studies of 
molecular motor uses periodic boundary conditions 
at $x_1$ and $x_2$.  Since there is no one-way permeable 
membrane, the boundary condition is $f^{ss}(x_1)=f^{ss}(x_2)\neq 0$. 
The cycle flux (i.e., mean velocity for a single motor) then is
\[
	J^{cycle} = \frac{e^{-\Psi(x_2)/\epsilon}
		- e^{-\Psi(x_1)/\epsilon}}
		{T_{x_1\rightarrow x_2}e^{-\Psi(x_2)/\epsilon}
		+T_{x_2\rightarrow x_1}e^{-\Psi(x_1)/\epsilon}}.
\]
The renewal process is then replaced by a semi-Markov 
process which can go both clockwise and counter-clockwise 
on a circle \cite{qian_xie_pre}.  Birth-death
processes on a circle will be the subject of a forthcoming
paper.

\subsection{Higher-order Kramers-Moyal expansion and 
Keizer's paradox}
\label{k_para}

Keizer's paradox was originally introduced to understand
a discrepancy between the infinitely long 
time behavior of a CME and its deterministic counterpart 
in terms of nonlinear ordinary differential 
equations (ODE) \cite{vellela_qian_07,keener_2}.  
The resolution is in the vast separation of time scales: 
the infinitely long time in the ODE is still a very short
time in the stochastic dynamics of the CME which involves 
``uphill climbing'' and ``barrier crossing''.  The same 
result also explains the discrepancy between the 
stationary distribution of a CME and the stationary 
distribution of the corresponding diffusion approximation
via a Fokker-Planck equation \cite{hanggi_84,vellela_qian_09}.  
This is now a well-understood subject, intimately 
related to the finite-time condition required in Kurtz's 
convergence theorem \cite{kurtz}: The convergence when
$V\rightarrow\infty$ is not uniform with respect to $t$.

We now offer a different, more explicit 
approach to illustrate the diffusion approximation problem.  The 
diffusion approximation of a CME in terms of 
a Fokker-Planck equation is actually a truncated 
Kramers-Moyal expansion up to 
$V^{-1}$ \cite{vanKampen,gardiner}.  
Naturally one can investigate the consequence of keeping 
the $V^{-2}$ term in the expansion:
\begin{equation}
	\frac{d}{dx}\left[
	\epsilon^2 a(x) \frac{d^2 u}{dx^2}
	+ \epsilon D(x)\frac{du}{dx} - b(x) u
		\right] = 0,
\label{2nd_order_km}
\end{equation}
where $\epsilon=1/V$.
Applying no-flux boundary condition at $x=\infty$,
we have
\begin{equation}
	\epsilon^2a(x) \frac{d^2 u}{dx^2}
	+ \epsilon D(x)\frac{du}{dx} - b(x) u = 0.
\label{bob_1}
\end{equation}
We apply the WKB method \cite{murray,omalley} by 
assuming the solution to Eq. (\ref{bob_1}) of the 
form
\begin{equation}
	u(x) = \exp\left[-\frac{1}{\epsilon}\phi_0(x)
			+ \phi_1(x) + \epsilon \phi_2(x) +
			\cdots \right].
\end{equation}
We then substitute the $u(x)$ into Eq. (\ref{bob_1}) and 
collect terms with the leading order $\epsilon^0$ to yield
\begin{equation}
     a(x)\left(\phi_0'(x)\right)^2 - 
	D(x)\phi_0'(x)-
	b(x) = 0.
\label{eq_4_wkb3}
\end{equation}
We note that if $a(x)\equiv 0$, then $\phi_0(x) = 
-\int_0^x (b(z)/D(z))dz$, as given in Eq. (\ref{the_Psi}). 
In this case, a root of $b(x)=0$, $x^*$, has
$\phi'_0(x^*)=0$ and $\phi''_0(x^*) = -b'(x^*)/D(x^*)$.
Hence, a stable fixed point of the ODE $dx/dt=b(x)$
corresponds to a local minimum of $\phi_0(x)$ and
a peak in the distribution $e^{-\phi_0(x)/\epsilon}$.

When $a(x)\neq 0$, Eq. (\ref{eq_4_wkb3}) still 
indicates that at $x^*$, the root of $b(x)=0$, $\phi'_0(x^*)=0$. 
Furthermore, differentiating Eq. (\ref{eq_4_wkb3}) with respect to
$x$ once, we obtain $\phi_0''(x^*)=-b'(x^*)/D(x)$. Therefore, 
the local behavior of $\phi_0(x)$ near a fixed point $x^*$ is 
independent of the higher order terms!  However, if 
$a(x)\neq 0$, then the global behavior of the solution to 
Eq. (\ref{eq_4_wkb3}) will have a non-negligible difference
from $-\int_0^x (b(z)/D(z))dz$.  This difference contibutes
to the difference $|\phi_0(x_1^*)-\phi_0(x_2^*)|$.

	Eq. (\ref{eq_4_wkb3}) is in fact a 
3rd-order truncated version of the exact equation for 
$\phi_0'(x)$ given by G. Hu \cite{hugang}:
\begin{equation}
	\mu_0(x)\left(e^{\phi_0'(x)}-1\right)
	+\lambda_0(x)\left(e^{-\phi_0'(x)}-1\right) = 0,
\label{exact}
\end{equation}
with $b(x)=\mu_0(x)-\lambda_0(x)$, $D(x)=(\mu_0(x)+\lambda_0(x))/2$, 
and $a(x)=(\lambda_0(x)-\mu_0(x))/6$.  The exact, non-trivial, solution to Eq. (\ref{exact}) is $\phi'_0(x)=\ln(\lambda_0(x)/\mu_0(x))$,
which is given in Eq. (\ref{asymp_4_pst}a) below.

\section{One-dimensional Birth-Death Processes: Stationary 
Distribution and Mean First Passage Time}

	We shall now be interested in the Kolmogorov forward 
equation (\ref{1dbdp}) for the 1-dimensional DGP. 
To be consistent with the macroscopic Law of Mass Action, 
we shall further assume both birth and death rates 
have asymptotic expansions in the limit of 
$V, n\rightarrow\infty$, $n/V\rightarrow z$:
\begin{eqnarray}
	V^{-1}u_{zV}(V) &=& \mu_0(z) + \frac{\mu_1(z)}{V} 
		+ \frac{\mu_2(z)}{V^2} + O\left(V^{-3}\right),
\\	
 	V^{-1}w_{zV}(V) &=& \lambda_0(z) + \frac{\lambda_1(z)}{V} 
		+\frac{\lambda_2(z)}{V^2} + 
		+ O\left(V^{-3}\right).
\end{eqnarray}

\subsection{Stationary distribution and local behavior near
a fixed point}

In terms of the birth and death rates $u_n(V)$ and $w_n(V)$,
the stationary distribution to Eq. (\ref{1dbdp}) is
\begin{equation}
	p^{st}_n = p_0^{st} \prod_{\ell=0}^{n-1}
                \frac{u_{\ell}(V)}{w_{\ell+1}(V)}.
\label{pstn}
\end{equation}
With large $V$, applying the lemma in Sec. \ref{math_lemma}, 
the asymptotic expansion for $p^{st}_n(V)$ is
\begin{eqnarray}
	\ln p^{st}_{xV}(V) &\approx& \ln f^{st}(x) 
		= V\int_0^x \ln\left(
		\frac{\mu_0(z)}{\lambda_0(z)}\right) dz  
\label{asymp_4_pst}
\\
	&-& \int_0^x \left(\frac{\lambda_1(z)}{\lambda_0(z)}
		- \frac{\mu_1(z)}{\mu_0(z)}\right) dz
		- \frac{\ln\left(\mu_0(x)\lambda_0(x)\right)}{2}
			+ O\left(V^{-1}\right),
\nonumber
\end{eqnarray}
in which we have neglected an $x$-independent 
term $\ln p^{st}_{0}(V)$.

	We shall point out that 
if one applies a {\em diffusion approximation} to the 
master equation (\ref{1dbdp}), we obtain the
diffusion equation (\ref{1dd}) with 
\begin{subequations}
\label{app_D_b}
\begin{eqnarray}
		D(x) &=& \frac{\mu_0(x)+\lambda_0(x)}{2V}+
		\frac{\lambda_1+\mu_1+\lambda'_0-\mu'_0}{2V^2},
		\ \textrm{ and } 
\\
		b(x) &=& \mu_0(x) - \lambda_0(x) + \frac{1}{V}
			\left(\mu_1-\lambda_1-\frac12(\lambda'_0+\mu'_0)\right).
\end{eqnarray}
\end{subequations}
Then the stationary distribution to the approximated
diffusion process is
\begin{eqnarray}
\label{diff_app}
	\ln \widetilde{f}(x) &=& 2V
		\int_0^x \frac{\mu_0(z)-\lambda_0(z)}
			{\mu_0(z)+\lambda_0(z)}\ dz 
\\
	&+& \int_0^x \frac{\lambda_0(4\mu_1+\lambda_0'-3\mu_0')
		-\mu_0(4\lambda_1+3\lambda_0'+3\mu_0')}
		{\left(\mu_0+\lambda_0\right)^2}\ dz + O\left(V^{-1}\right),
\nonumber	
\end{eqnarray}
which is different from Eq. (\ref{asymp_4_pst}),  
even in the leading order.

	However, it is easy to verify that the leading-order
terms in the indefinite integrals in (\ref{asymp_4_pst}) 
and (\ref{diff_app}), as functions of $x$, have matched 
locations for their extrema as well as their curvatures at 
each extrema.  This is because both
\begin{equation}
	\frac{d\ln f^{st}(x)}{dx} = V\ln\frac{\mu_0(x)}{\lambda_0(x)}
	\  \textrm{ and } \ 	
	\frac{d\ln\widetilde{f}(x)}{dx} = 
		2V\frac{\mu_0(x)-\lambda_0(x)}{\mu_0(x)+\lambda_0(x)}
\end{equation}
are zero at the root of $b(x)=\mu_0(x)-\lambda_0(x)$. 
One can further check that both have identical slopes at their 
corresponding zeros.

	Therefore, near a stable fixed point $x^*$ of 
$dx/dt=\mu_0(x)-\lambda_0(x)$: $\mu_0(x^*)=\lambda_0(x^*)$ 
and $\mu'(x^*)<\lambda'(x^*)$, both approaches 
yield a same Gaussian process with diffusion equation 
\begin{equation}
   \frac{\partial f(\xi,t)}{\partial t}
	= \frac{\lambda_0(x^*)}{V}
	\frac{\partial^2 f(\xi,t)}{\partial\xi^2}
	-  \frac{\partial}{\partial \xi}
	\left((\mu'(x^*)-\lambda'(x^*)) \xi f(\xi,t)\right).
\end{equation}
where $\xi=x-x^*$ is widely called {\em fluctuations} in 
statistical physics. This is Onsager-Machlup's Gaussian 
fluctuation theory in the linear regime
\cite{onsager,keizer_jmp,qian_prsa_01,qian_review_3}.

\subsection{Mean first passage time and diffusion approximation}
\label{sec_kbe}

	Corresponding to the Kolmogorov forward equation in
Eq. (\ref{1dbdp}), the Kolmogorov backward equation for the 
birth-death process is
\begin{equation}
	\frac{dg_n}{dt} = w_ng_{n-1} - (u_n+w_n)
			g_n + u_ng_{n+1}.
\label{kbe_4_bdp}
\end{equation}
Then, $T_n$ ($0\le n\le n_2$), the mean first passage time
(MFPT)  arriving at $n_2$,
starting at $n$ with a reflecting boundary at
$0$,  satisfies the 
inhomogeneous equation
\begin{equation}
	w_nT_{n-1} - (u_n+w_n)
			T_n + u_nT_{n+1} = -1,
\end{equation}
with the boundary conditions
\begin{equation}
	T_0=T_{-1} \ \textrm{ and } \ T_{n_2} = 0.
\end{equation}
The solution can be found in many places, {\it e.g.},
Ch. XII in \cite{vanKampen} and Equn. 31 in \cite{doering_1}.  
The result is most compact when expressed 
in terms of the stationary distribution
$p^{st}_n(V)$ in Eq. (\ref{pstn}):
\begin{equation}
 T_{n\rightarrow n_2} =
        \sum_{m=n+1}^{n_2}\sum_{\ell=0}^{m-1}\frac{p^{st}_{\ell}(V)}
                                {w_m(V)p^{st}_m(V)}.
\label{T12_exact}
\end{equation}
In the limit of large $V$, one has the asymptotic 
expression for $T_{n_1\rightarrow n_2}$ 
(see Sec. \ref{sec.vii_b}):
\begin{equation}
   V\int_{x}^{x_2}
	\frac{\ln\left[\frac{\lambda_0(z)}{\mu_0(z)}\right]}
	{\left[\frac{\lambda_0(z)}{\mu_0(z)}-1\right]}
	 \frac{e^{V\Phi(z,V)}}{\lambda_0(z)} dz\int_{0}^{z}
   \frac{\ln\left[\frac{\mu_0(y)}{\lambda_0(y)}\right]}
	{\left[\frac{\mu_0(y)}{\lambda_0(y)}-1\right]}
	 e^{-V\Phi(y,V)} dy,
\label{asym_4_T}
\end{equation}
in which $x=n/V$, $x_2=n_2/V$, and
\begin{equation}
	\Phi(x,V) = -\frac{1}{V}\ln p^{st}_{xV}(V)
           = \phi_0(x)+\frac{1}{V}\phi_1(x)
			+ O\left(V^{-2}\right),
\end{equation}
given in Eq. (\ref{asymp_4_pst}). 

	Comparing Eqs. (\ref{asym_4_T}) and (\ref{mfpt}), we 
see that the effective ``potential function''
\begin{eqnarray}
	\widetilde{\Psi}(x,V) &=& \Phi(x,V)
		+\frac{1}{V}\ln\frac{\mu_0(x)/\lambda_0(x)-1}
			{\ln\mu_0(x)-\ln\lambda_0(x)},
\label{tilde_Psi}
\end{eqnarray}
and effective diffusion coefficient
\begin{equation}
	\widetilde{D}(x) = \frac{1}{\mu_0(x)}\left(\frac{\mu_0(x)-\lambda_0(x)}
		{\ln\mu_0(x)-\ln\lambda_0(x)}\right)^2.
\label{tilde_D}
\end{equation}
We note that near $\mu_0(x)=\lambda_0(x)$, 
\begin{equation}
	\widetilde{D}(x) \approx \lambda_0(x)\left[1+\frac{(\lambda_0-\mu_0)^2}
		{12\lambda_0^2}\right].
\label{eq_30}
\end{equation}
Following Eq. (\ref{the_Psi}), 
Eqs. (\ref{tilde_Psi}) and (\ref{tilde_D}) imply that 
\begin{equation}
	\widetilde{b}(x) = \left(\frac{1-\lambda_0(x)/\mu_0(x)}
			{\ln\mu_0(x)-\ln\lambda_0(x)}\right)
			\left(\mu_0(x)-\lambda_0(x)\right).
\end{equation}
As we shall show, even though the general formula for the MFPT
given in Eq. (\ref{asym_4_T}) has a complex expression for
the effective diffusion coefficient $\widetilde{D}(x)$ and
effective drift $\widetilde{b}(x)$, the Kramers-like 
formular for barrier crossing, Eq. (\ref{new_k_eq}) which 
only involves $D(x)$ at the peak of the potential function, is 
simple and recognizable.  

Noting the disagreement between the correct 
asymptotic Eq. (\ref{asymp_4_pst}) and the stationary
distribution (\ref{diff_app}) from the diffusion approximation 
with $D(x)$ and $b(x)$ given in Eq. (\ref{app_D_b}),
H\"{a}nggi {\em et al.} proposed an alternative
diffusion equation with
\begin{equation}
	D^{hgtt}(x) =	\frac{\mu_0(x)-\lambda_0(x)}{\ln\mu_0(x)-\ln\lambda_0(x)},
 \ \
	b(x) = \mu_0(x)-\lambda_0(x),
\label{eq_31}
\end{equation}
as a more appropriate approximation
for the 1-dimensional Eq. (\ref{1dbdp}) \cite{hanggi_84}.
While the H\"{a}nggi-Grabert-Talkner-Thomas diffusion yields
the correct leading order stationary distribution
(\ref{asymp_4_pst}), we note that the $D^{hgtt}(x)$ 
is different from the $\widetilde{D}(x)$ in Eq. (\ref{tilde_D}).  
In fact,
\begin{equation}
	\frac{D^{hgtt}(x)}{\widetilde{D}(x)} = 
		\frac{\ln\mu_0-\ln\lambda_0(x)}
			{1-\lambda_0/\mu_0} < 1 \textrm{ when }
		\frac{\lambda_0}{\mu_0} > 1; \ \ > 1 \textrm{ when } 
		\frac{\lambda_0}{\mu_0} < 1,
\end{equation}
and near $\mu_0(x)=\lambda_0(x)$, 
\begin{equation}
  D^{hgtt}(x) = \lambda_0\left[1+\frac{\mu_0-\lambda_0}{2\lambda_0}
		-\frac{(\lambda_0-\mu_0)^2}
		{12\lambda_0^2}\right].
\label{eq_33}
\end{equation}

	The H\"{a}nggi-Grabert-Talkner-Thomas 
diffusion process with diffusion and drift given 
in Eq. (\ref{eq_31}) is the only diffusion process 
that yields the correct deterministic limit 
$b(x)$ and asymptotically correct stationary
distribution (\ref{asymp_4_pst}).  However, the 
discrepancy between $D^{hgtt}(x)$ and $\widetilde{D}(x)$
(Eqs. \ref{eq_30} and \ref{eq_33})
indicates that it will not, in general, give the
asymptotically correct MFPT.  Therefore, diffusion 
processes have difficulties in approximating {\em both} the 
correct stationary distribution {\em and} the 
correct MFPT of a birth-death 
process in the asymptotic limit of large $V$.  This 
conclusion has been dubbed {\em diffusion's dilemma} \cite{qian_review_3,zhou_qian_pre}.

\subsection{Kramers' formula and MFPT for barrier crossing}

	With a correct potential fuunction $\phi_0(x)$,
the difference between the discrete DGP with asymptotic 
large $V$ and continuous approximation disappears in 
the computations of MFPT for Kramers problem, i.e., 
barrier crossing.  This is because, as we have 
demonstrated, all different approximations can 
preserve local curvatures at the stable and unstable
fixed points.  At a fixed point $x^*$ of $b(x)$: $b(x^*)=0$, 
and the $D(x^*)=\lambda_0(x^*)=\mu_0(x^*)$ in 
both equations (\ref{eq_30}) and (\ref{eq_33}).
And with a correct ``barrier height'' and local
curvatures at the extrema, the Kramers formula is completely determined. 

	In fact, applying Laplace's method and considering
an energy barrier between $x_1^*=n_1^*/V$ and $x_2^*=n_2^*/V$,
located at $x^{\ddag}$, Eq. (\ref{asym_4_T}) can be simplied 
into (see Sec. \ref{sec.vii_c} for details)
\begin{equation}
	T^* = T_{n_1\rightarrow n_2} = 
		\frac{2\pi}{\lambda_0(x^{\ddag})
		\sqrt{\phi_0^{''}(x_1^*)|\phi_0^{''}(x^{\ddag})|}} e^{V[\Phi(x^{\ddag},V)-\Phi(x_1^*,V)]} 
		\left\{ 1 + O\left(\frac{1}{V}\right) \right\},
\label{new_k_eq}
\end{equation} 
in which $\Phi(x,V)=\phi_0(x)+(1/V)\phi_1(x)$.
Note that Eq. (\ref{new_k_eq}) contains a $(1/V)\phi_1(x)$
term in exponent.  This is a key result of the present
paper, a new feature for the DGP.

	Note also that the ambiguity discovered in \cite{doering_2}
is associated with MFPT with both starting and end points 
within a same basin of attraction.  It does not appear 
in Kramers' formula for inter-basin transition. 

	The $T^*$ given in Eq. (\ref{new_k_eq}) 
plays an all-important role in the dynamics 
with mutiple attractors.  It divides the local,
intra-basin dynamics from the stochastic,
inter-basin jump process.
The barrier 
\begin{equation}
	V\left(\Phi(x^{\ddag},V)-\Phi(x^*_1,V)\right)=
  V\{\phi_0(x^{\ddag})-\phi_0(x^*_1)\}+
 \{\phi_1(x^{\ddag})-\phi_1(x^*_1)\},
\end{equation}
thus the time $T^*$, can be increasing or 
decreasing with $V$, depending on
$\phi_0(x^{\ddag})-\phi_0(x^*_1)>0$ or $<0$.
This distinction leads to the concept of {\em nonlinear
bistability} vs. {\em stochastic bistability}
\cite{bishop_qian_bj_10,qian_review_2,zhang_ge_qian}.

\section{Nonlinear Bifurcation and Stochastic Phase 
Transition: A Potential Function Perspective}

{\em Bifurcation} is one of the most important characteristics 
of nonlinear dynamical systems \cite{perko,kuznetsov}.  
Therefore, a nonlinear stochastic dynamic theory cannot be 
complete without a discussion of its stochastic counterpart.    
This is still a developing area in random dynamical
systems \cite{arnold} and stochastic processes \cite{schuss}.  
In fact, the very notion of stochastic bifurcation has 
at least two definitions: the P- (phenomenological) and 
D- (dynamic) bifurcations \cite{arnold}.  We shall not 
discuss the fundamental issues here; rather
provide some observations based on applied mathematical
intuition.  This discussion is more consistent 
with the P-bifurcation advocated by 
E.C. Zeeman \cite{zeeman}.

The canonical bifurcations in one-dimensional nonlinear 
dynamical systems are {\em transcritical},
{\em saddle-node}, and {\em pitchfork} bifurcations
\cite{perko}.  The saddle-node bifurcation and its
corresponding stochastic model have been extensively
studied in terms of the Maxwell construction
\cite{ge_qian_jrsi}. The key is to realize
the separation of time scales, and the different orders
in taking limits $V\rightarrow\infty$ and $t\rightarrow\infty$.  
The steady states
of deterministic nonlinear ODEs are initial
value dependent; but the steady state distribution
of the stochastic counterpart, in the limit of
$V\rightarrow\infty$, is unique and independent
of the initial value.  An ODE finds a ``local minimum''
of $\Psi(x)$ in the infinite time while its stochastic 
counterpart finds the ``global minimum'' at its infinite time.

In this section, We shall mainly discuss the transcritical 
bifurcation which has not attracted much attention
in the past.  We shall show that in certain cases, it is in 
fact intimately related to the extinction phenomenon and 
 Keizer's paradox.  

\subsection{Transcritical bifurcation}

	The normal form of transcritical 
bifurcation is \cite{perko}
\begin{equation}
	\dot{x} = b(x) \approx \mu (x-x^*) - (x-x^*)^2 
			\ \textrm{ near $x=x^*>0$},
\label{trans_cri_b}
\end{equation}
in which the locale of bifurcation is at $x=x^*$; It 
occurs when $\mu=0$.  Transcritical bifurcation
is a local phenomenon and the far right-hand-side
of (\ref{trans_cri_b}) is the Taylor expansion 
of $b(x)$ in the neighbourhood of $x=x^*$.   Let us 
assume that the system's lower bound is $0$; for
an ODE to be  meaningful to population dynamics, 
the $b(0)$ has to be non-negative.  This implies
$b(x)$ has another, stable fixed point $x_1^*$, 
$0 \le x_1^* < x^*$ when $|\mu|$ is sufficiently small.

At the critical value of $\mu=0$, the 
steady state at $x=x^*$ has the form of
$\dot{x}=-(x-x^*)^2$.  Hence the corresponding
potential function, near $x=x^*$, will be
\begin{equation}
	\Psi(x;\mu=0) = -\int_{x^*}^x \frac{b(z)}{D(z)}dz
		  = \frac{(x-x^*)^3}{3D(0)} + O\left((x-x^*)^4\right).
\end{equation}
Since $D(x)>0$, the $\Psi(x)$ is neither a minimum
nor a maximum at $x=x^*$.  There is a minimum of $\Psi(x)$ 
at $x_1^*$.  Now for sufficiently small $\mu\neq 0$, 
a pair of minimum and maximum develop in the neighbourhood
of $x^*$, approximately at $x^*$ and $x^*+\mu$.
Then by continuity, the newly developed minimum 
of $\Psi(x,\mu)$ must {\em not} be lower than $\Psi(x_1^*,\mu)$. 
In other words, the stationary distribution of 
the stochastic dynamics, in the limit of $V\rightarrow\infty$,
will not be in the neighbourhood of the location of a
transcritical bifurcation.

\begin{figure}[t]
\centerline{\includegraphics[width=3.5in]{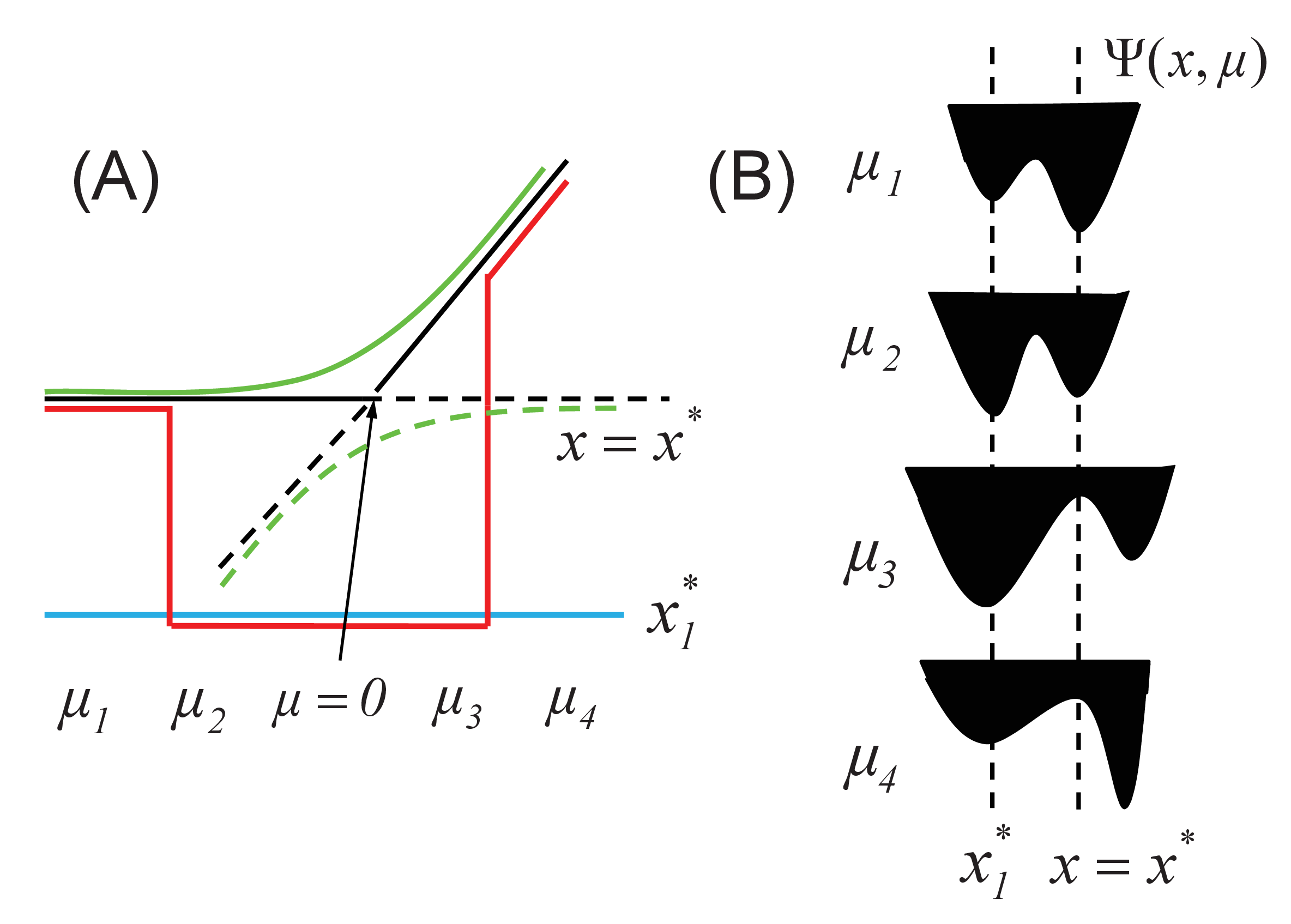}}
\caption{Transcritical bifurcation ocurrs between 
$\mu_2$ and $\mu_3$ at $x=x^*$, as shown by the solid 
and dashed black lines in (A), and corresponding 
stochastic potential $\Phi(x)$ shown in (B).  
Phase transition(s) can occur between $\mu_1$ 
and $\mu_2$, and between $\mu_3$ 
and $\mu_4$, when the global minimum of the potential 
function $\Psi(x,\mu)$ switches between at $x_1^*$ (blue) 
and near $x^*$ (solid black).  The global minimun is
represented by the red curve.   The phase transition,
therefore, is not really associated with the
transcritical bifurcation.  A very minor ``imperfection'' 
can lead the two black lines, solid and dashed, to 
become the two green lines.  While the transcritical 
bifurcation disappeared, the phase transitions are still 
present.  The latter phenomenon is structurally stable
while the former is not \cite{zeeman}. The bifurcation
is local while the phase transitions are global. 
}
\label{fig_2}
\end{figure}

	However, as in the case of saddle-node bifurcation,
phase transition might occurs for larger $|\mu|$.  Note
that the local minimum associated with the transcritical 
bifurcation could become the global minimum of $\Phi(x,\mu)$.
Then that occurs, there is a phase transition.
This is illustrated in Fig. \ref{fig_2}.  Note however, 
the phase transitions are not associated with the 
transcritical bifurcation {\em per se}: It is really the 
competition between two stable fixed points, the solid green 
line and the blue line, $x_1^*$.

	If the $x=x^*$ happens at the boundary of the domain
of $x$, we have a more interesting scenario.
Consider the birth-death system with
$u_n(V) = k_1n $ and $w_n(V)= k_{-1}n(n-1)/V+k_2n$.
Then  the corresponding $\mu_0(x)=k_1x$,
$\lambda_0(x)=k_{-1}x^2+k_2x$, and 
\begin{equation}
	b(x) = \mu_0(x)-\lambda_0(x) = (k_1-k_2)x-k_{-1}x^2,
	\ \  x\ge 0.
\end{equation} 
The ODE $\dot{x}=b(x)$ has a transcritical bifurcation at 
$x=0$ when $k_1-k_2=0$.  When $k_1<k_2$, the
system has only a stable fixed point at $x=0$;
when $k_1>k_2$, it has a stable fixed point at 
$x=\frac{k_1-k_2}{k_{-1}}$, and $x=0$ is a unstable
fixed point.

	However, the stochastic stationary distribution
has a probability 1 at $n=0$, i.e., extinction, for 
any value of $k_1-k_2$.  This is Keizer's paradox 
\cite{vellela_qian_07}.

	Therefore, when a transcritical bifurcation
ocurrs at the boundary of a domain, the stochastic 
steady state exhibit no discontinuous ``phase transition''.
Rather, the boundary is an absorbing state.  On the
other hand, when a transcritical bifurcation ocurrs in
the interior of the domain, there might not be phase
transition associated with it.  Transcritical
bifurcation and phase transition are two different 
phenomena.

\begin{figure}[t]
\centerline{\includegraphics[width=2.75in,angle=-90]{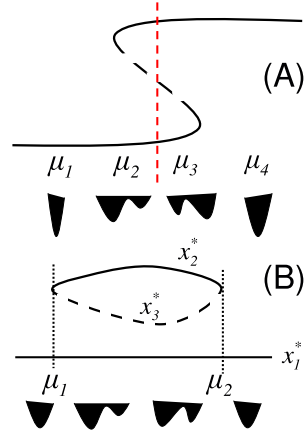}}
\caption{(A) Two saddle-node bifurcations together give
the catastrophe phenomenon.  At the four different parameter
$\mu$ values, the stochastic potential has a single 
minimum ($\mu_1$), then passing through saddle-node 
bifurcation to have two minima, with the lower one being
the global minimum ($\mu_2$).  Then at $\mu_3$, the 
global minimum is the upper one.  And finally at
$\mu_4$, saddle-node bifrucation leads again to a single steady
state.  The stochastic phase transition occurs between
$\mu_2$ and $\mu_3$, denoted by the red dashed line known
as the Maxwell construction \cite{ge_qian_jrsi}.  The
stochastic potential $\Psi(x)$ corresponding to the four
$\mu$'s are illustrated below the bifurcation diagram.
(B)  Two saddle-node bifurcations
occur at $\mu_1$ and $\mu_2$.  In this case, however, 
there is no phase transition, as illustrated by
the given $\Psi(x)$ below the bifurcation diagram.
}
\label{fig_3}
\end{figure}

\subsection{Saddle-node bifurcation}

	The normal form of saddle-node bifurcation 
is \cite{perko}
\begin{equation}
	\dot{x} = b(x) \approx \mu - x^2 
			\ \textrm{ near $x=0$},
\label{sa_no_b}
\end{equation}
in which the locale of bifurcation is again at $x=0$ 
and it occurs when $\mu=0$. 
Again, it is clear that the stochastic phase transition
is not associated with a single saddle-node bifurcation
event {\em per se}.  However, it is necessitated by the {\em two} 
saddle-node bifurcation events in a catastrophe 
phenomenon, which has a deeper topological root.  This is 
illustrated in Fig. \ref{fig_3}.

\section{\blue{van't Hoff-Arrhenius Analysis}}

	There is a deep connection between the potential 
function $\Psi(x,V)$ and the theory of 
thermodynamics.  
In this section we provide a brief discussion of the
subject, which is yet to be fully developed.
 
	In classical thermodynamics, there
is a decomposition of ``free energy into 
enthalpy and entropy''.  The canonical definitions are
\begin{eqnarray}
	\widetilde{\phi}_0(x,V) &=& \left(\frac{\partial 
			\left(V\Phi(x,V)\right)}
				{\partial V}\right)_x,
\label{vanthoff}
\\
        \widetilde{\phi}_1(x,V) &=& \left(\frac{\partial \Phi(x,V)}
				{\partial (1/V)}\right)_x.
\label{entropy_def}
\end{eqnarray}
It is easy to verfiy that $\widetilde{\phi}_0+(1/V)\widetilde{\phi}_1
=\Phi(x,V)$.  Noting an analogue between temperature $T$ and 
$1/V$, Eq. (\ref{vanthoff}) is known as the van't Hoff equation for 
enthalpy in classic thermodynamics, and Eq. (\ref{entropy_def}) 
is the definition of entropy: $S=-\left(\partial G/\partial T\right)_{P,N}$.
Then $\widetilde{\phi}_0$ and $-\widetilde{\phi}_1$ are the 
``enthalpic'' and ``entropic'' components of ``free energy'' 
$\Phi$.  Knowing one of $\Phi(x)$, $\phi_0(x)$, and 
$\phi_1(x)$, one can determine the other two 
functions within an additive constant: In classical
thermodynamics, they are all called ``thermodynamic
potentials'', and they are equivalent.

	Since in the limit of $V\rightarrow\infty$, 
$\Phi(x,V)=\phi_0(x)+(1/V)\phi_1(x)+ o\left(V^{-1}\right)$,
we have
\begin{equation}
\lim_{V\rightarrow\infty} \widetilde{\phi}_0(x,V) = \phi_0(x),
\ \textrm{ and } \ 
\lim_{V\rightarrow\infty} \widetilde{\phi}_1(x,V) = \phi_1(x).
\end{equation}
Furthermore, for $\Phi(x,V)$ being a continuous function of
$V$, if $\widetilde{\phi}_0(x,V)=\phi_0(x)$ is
independent of $V$, then
$\widetilde{\phi}_1(x,V)=\phi_1(x)$ is independent of $V$.
This is because 
\begin{equation}
	\frac{\partial}{\partial V}\ \widetilde{\phi}_1(x,V)
	= -V\frac{\partial}{\partial V}\ \widetilde{\phi}_0(x,V).
\end{equation}
In this case, $\Phi(x,V)$ is a linear function of $1/V$
with slope $\phi_1(x)$ and intersect $\phi_0(x)$.
Systems with $V$ independent $\widetilde{\phi}_0(x,V)$
are ``thermodynamically'' simple systems.

\subsection{\blue{Decomposing the potential \boldmath{$\Phi(x,V)$} into \boldmath{$\widetilde{\phi}_0(x,V)$}\\
\boldmath{$+(1/V)\widetilde{\phi}_1(x,V)$}}}
	 
The stochastic potential for a one-dimensional birth-death
processes has the explicit expression:
\begin{equation}
	\Phi(x,V) = -\frac{1}{V}\ln p^{st}_{xV}(V)
		  = \frac{1}{V}
		\sum_{\ell=0}^{xV-1}\ln\frac{w_{\ell+1}(V)}{u_{\ell}(V)},
\label{eq_43}
\end{equation}
in which we have neglected, on the right-hand-side, an added term 
$-(1/V)\ln p^{st}_0(V)$ which is a function of only $V$, not $x$.  
Then following Eqs. (\ref{vanthoff}) and (\ref{entropy_def}), we 
have
\begin{eqnarray}
	\widetilde{\phi}_0(x,V) &=& -\left[\left(\frac{\partial\ln p^{st}_n(V)}
		{\partial V}\right)_n+
		x\left(\frac{\partial\ln p^{st}_n(V)}{\partial n}\right)_V
		\right]_{n=xV},
\\
	\widetilde{\phi}_1(x,V) &=& \left[V\left(\frac{\partial\ln p^{st}_n(V)}
		{\partial V}\right)_n+
		Vx\left(\frac{\partial\ln p^{st}_n(V)}{\partial n}\right)_V
		-\ln p^{st}_n(V)\right]_{n=xV},
\end{eqnarray}
in which we shall substitute the derivative with respect to $n$ by
a difference:
\begin{equation}
	\left(\frac{\partial\ln p^{st}_n(V)}{\partial n}\right)_V
	    = \ln p^{st}_n(V) - \ln p^{st}_{n-1}(V)
	    = \ln\frac{u_{n-1}}{w_n}.
\end{equation}
Thus, the decomposition:
\begin{eqnarray}
	\widetilde{\phi}_0(x,V) &=& -\left[\sum_{\ell=0}^{n-1}
		\frac{d}{dV}\left(\ln\frac{w_{\ell+1}(V)}{u_{\ell}(V)}\right)
		+x\ln\frac{u_{n-1}(V)}{w_n(V)}
		\right]_{n=xV},
\label{eq_47}
\\
	\widetilde{\phi}_1(x,V) &=& \left[\sum_{\ell=0}^{n-1}
		\frac{d}{d\ln V}\left(\ln\frac{w_{\ell+1}(V)}
		{u_{\ell}(V)} \right)
		+Vx\ln\frac{u_{n-1}(V)}{w_n(V)}
		-\ln p^{st}_n(V)\right]_{n=xV}.
\end{eqnarray}

\subsection{\blue{Chemical reaction systems with a linear 
\boldmath{$1/V$} dependent potential}}

	What type of CMEs will have a simple, $V$
independent $\widetilde{\phi}_0(x,V)$ and
$\widetilde{\phi}_1(x,V)$?  Let us consider a 
simple example:  the Poisson distribution with 
$u_{\ell}(V)=\alpha V$, $w_{\ell}(V)=\beta\ell$, and
\begin{equation}
	p_n(V) = \frac{\left(x^*V\right)^n}{n!}e^{-x^*V}, \ \ 
		x^*=\frac{\alpha}{\beta}.
\end{equation}
We then have
\begin{eqnarray}
	\widetilde{\phi}_0(x,V) &=& x\ln\frac{x}{x^*}-x+x^*,
\\
	\widetilde{\phi}_1(x,V) &=& \left(\ln n! -n\ln n + n\right)_{n=xV}
		\approx \ln\sqrt{2\pi x} + \textrm{const.}
\end{eqnarray}
Therefore, the Poisson distribution has a 
stochastic potential with linear $1/V$ dependence.

	On the other hand, the binomial distribution with
$u_{\ell}(V)=k_+(e_tV-\ell)$, $w_{\ell}(V)=k_-\ell$, and
\begin{equation}
	p_n(V) = \frac{N!}{n!(N-n)!}\frac{\theta^n}{(1+\theta)^N},
	\ \ \theta = \frac{k_+}{k_-}, \ N=e_tV.
\end{equation}
Then,
\begin{eqnarray}
	\widetilde{\phi}_0(x,V) &=& \left[\sum_{\ell=0}^{n-1}\frac{e_t}{e_tV-\ell}
			-x\ln\frac{\theta(e_tV-n+1)}{n}\right]_{n=xV}
\nonumber\\
	&\approx& e_t\ln\frac{e_t-x}{e_t}-x\ln\frac{\theta(e_t-x)}{x}
			+\frac{x}{2V(e_t-x)},
\label{eq_53}
\\
	\widetilde{\phi}_1(x,V) &=& \left[-\sum_{\ell=0}^{n-1}\frac{e_tV}{e_tV-\ell}
			+xV\ln\frac{e_tV-n+1}{n}+\ln n!+\ln(e_tV-n)!\right]_{n=xV}
\nonumber\\
	&\approx& \ln\sqrt{(e_t-x)x} + \textrm{const.}
\end{eqnarray}
We see that in Eq. (\ref{eq_53}), there is a $(1/V)$
dependence.  Hence even the simple binomial distribution
has a $V$-dependent $\widetilde{\phi}_0(x,V)$.  Note
that because the potential function $\Phi(x,V)$ defined
in Eq. (\ref{eq_43}) is not a continuous function of 
$V$, $V$-independent $\widetilde{\phi}_1(x,V)$ no
longer dictates $V$-independent $\widetilde{\phi}_0(x,V)$.

\section{\blue{Discussion}}

	Multi-dimensional birth-death processes have
become in recent years a fundamental theory for stochastic,
nonlinear dynamical systems with individual, ``agent-based''
stochastic nonlinear behavior, and emergent 
long-time discontinuous stochastic evolution \cite{qian_review_2,
qian_review_3}.  The dynamics of such a process is
intimately related to both nonlinear ordinary
differential equations and multi-dimensional diffusion
processes.  In the present paper, we developed a systematic
study for the simplest, one-dimensional system.  Each of
the equations (\ref{the_ode}), (\ref{1dd}) and (\ref{1dbdp})
plays a role in the theory of a birth-death 
process.  

	One of the main recent applications of this type of 
dynamic models is in cellular biochemistry in a mesoscopic
volume, the Delbr\"{u}ck-Gillespie process (DGP) whose
Kolmogorov forward equation is widely known as the chemical 
master equation (CME), and whose stochastic trajectories
can be sampled following the Gillespie algorithm.
One special feature of a DGP is a parameter $V$, the 
system size. When $V\rightarrow\infty$, the trajectory 
of a DGP becomes the solution to an ODE. 

The application of the CME and the birth-death processes
to unimolecular, linear reaction systems can be found
in \cite{othmer_05,qian_jcp_06,huisinga_07}.  Steijaert 
{\em et al.} \cite{steijaert_10} have also presented a 
coherent summary for the CME with a single variable.  
In particular, they studied the metastability in the 
thermodynamic limit ($V,n\rightarrow\infty$; $n/V=x$)
associated with the Maxwell construction and stochastic
phase transition \cite{ge_qian_prl,ge_qian_jrsi}.

In the present work, we first studied the asymptotic
stationary solution of the 1-d stochastic dynamics for large 
system size $V$, $e^{-V\phi_0(x)-\phi_1(x)}$.  We 
then expressed the asymptotic mean first passage time
(MFPT) for a 1-d DGP in terms of the stochastic potential
$\Phi(x,V)\approx\phi_0(x)+(1/V)\phi_1(x)$.  With 
these two results, we obtain an effective diffusion 
coefficient, Eq. (\ref{tilde_D}) and a potential. 
A diffusion process defined by these two is significantly 
different from the diffusion approximation proposed by
H\"{a}nggi {\em et al.} \cite{hanggi_84} (comparing
Eqs. \ref{tilde_D} with \ref{eq_33}).  The latter does
yield the correct stationary distribution.  Together, 
our analysis shows that no single diffusion process with 
continuous stochastic paths can be globally asymptotically
accurate for the birth-death processes in general.  
A situation we have called ``diffusion's dilemma''
\cite{qian_review_3,zhou_qian_pre}.

In the limit of large system size $V$, a DGP process
with multistability (or multimodality) exhibits very 
different dynamics on different time scale:  
For $t\ll$ a critical $T^*$,
the dynamics is a continuous Gaussian process whose
mean value follows a deterministic linear dynamic.  
We call this {\em intra-basin dynamics}.  For $t> T^*$,
however, an {\em inter-basin dynamics} constitutes
Markov jump process that moves from an attractor to
another attractor.   The MFPT gives an estimation for the
$T^*$, which is in fact the reciprocal of the absolute
value of the second largest eigenvalue of the system. 
The largest eigenvalue is always zero for a Markov process \cite{hanggi_84,vellela_qian_07,vellela_qian_09}. 

	There is a deep relation between the nonlinear
bifurcation phenomenon, which is when $t\rightarrow\infty$
in a deterministic ODE but still $t\ll T^*$, and 
stochastic phase transition, which is related to
$t\rightarrow\infty$ and $t\gg T^*$.  In the theory of
DGP, phase transition can be studied in terms of the 
stochastic potential $\Phi(x,V)$ in the limit of
$V\rightarrow\infty$.  The limit of $N\rightarrow\infty$
followed by $t\rightarrow\infty$ is widely called 
{\em quasi-steady state} which is initial value dependent; 
while the limit $t\rightarrow\infty$
followed by $N\rightarrow\infty$ is unique, except at the
critical point of a phase transition: Phase transition 
occurs in the stationary distribution of an infinitely 
large system.  In the present paper for 1-d systems,
we have analyzed both transcritical bifurcation and
saddle-node bifurcation.  They are local
behaviors while a phase transition is a global
phenomenon.  However, the catastrophe phenomenon
necessitates a phase transition in terms of the 
Maxwell construction \cite{ge_qian_prl} 
(see Fig. \ref{fig_3}a).

	Finally, we also suggested an interesting connection
between the mathematical theory of DGP and the classical
thermodynamics in terms of the van't Hoff-Arrhenius
decomposition of free energy into enthalpy and entropy.

\section{\blue{Mathematical Details}}

\subsection{A lemma}
\label{math_lemma}

	Repeatedly using the definition of Riemann integration
\begin{equation}
	\lim_{V\rightarrow\infty}
	\sum_{\ell=0}^{xV-1} \frac{1}{V}\ F\left(\frac{\ell}{V}\right)
		= \int_0^x F(z)dz,
\end{equation}
we have for a smooth function $F(x)$: 
\begin{eqnarray}
	\sum_{\ell=0}^{xV-1} \frac{1}{V}\ F\left(\frac{\ell}{V}\right)
	&=& \int_0^x F(z)dz +\sum_{\ell=1}^{\infty}
		\frac{\nu_{\ell}}{V^{\ell}}
		\int_0^x F^{(\ell)}(z)dz
\\
	&=& \int_0^x F(z)dz -\frac{F(x)-F(0)}{2V}
		+ \frac{F'(x)-F'(0)}{12V^2} + O\left(V^{-3}\right).
\nonumber
\end{eqnarray}
where the coefficients $-\frac{1}{2}$, $\frac{1}{12}$, $0$, $\cdots$,
are the coefficients of the Taylor expansion of
\begin{equation}
	\frac{x}{e^x-1}-1 = -\frac{x}{2}+\frac{x^2}{12}+\cdots
			= \sum_{\ell=1}^{\infty} \nu_{\ell}x^{\nu}.
\end{equation}

	Now consider $F(x)$ that is not only a function of $x$, but
also a fuction of $V$. For example $\ell(\ell-1)(\ell-2)/V^3=x(x-1/V)(x-2/V)=F(x,V)$. 
In the limit of $V\rightarrow\infty$, one has the asymptotic 
expansion
\begin{equation}
	F(x,V)= F_0(x) + \frac{1}{V}F_1(x)
		+ \frac{1}{V^2}F_2(x)+\cdots,
\end{equation}
then 
\begin{eqnarray}
	&& \sum_{\ell=0}^{xV-1} \frac{1}{V}\ F\left(\frac{\ell}{V},V\right)
	= \int_0^x F_0(z)dz + \frac{1}{V}
		\left(-\frac{F_0(x)-F_0(0)}{2}
		+ \int_0^x F_1(z)dz \right)
\nonumber\\[-6pt]
\\[-4pt]
	&&+ \frac{1}{V^2}\left(\frac{F_0'(x)-F_0'(0)}{12} 
		-\frac{F_1(x)-F_1(0)}{2}
		+\int_0^x F_2(z)dz \right) + O\left(V^{-3}\right).
\nonumber
\end{eqnarray}

\subsection{Detailed derivation for Eq. (\ref{asym_4_T})}
\label{sec.vii_b}

By Taylor expansion and Riemann summation, we
shall first show Eq. (\ref{the_result00}) below. 

We starts with
\[
F(V)=V^2\sum_{m=m_1}^{m_2}\sum_{n=0}^{m} f(m,n)
         \hspace{3in}
\]
\begin{equation}
 \times \int_{s(m)}^{s(m+1)}\int_{s(n)}^{s(n+1)}
    e^{V[g(m,n)+g_x(m,n)\Delta x+g_y(m,n)\Delta y]} \ dydx,
\label{def_FV}
\end{equation}
in which
$g_x(x,y)$ and $g_y(x,y)$ denote the partial derivatives
$\big(\partial g/\partial x\big)_y$ and 
$\big(\partial g/\partial y\big)_x$.  Furthermore,
$\Delta x_m=x-s(m)$, $\Delta y_n=y-s(n)$, and 
$s(m) =  a+\frac{m-m_1}{m_2-m_1}(b-a)$ with
$(b-a)=V^{-1}(m_2-m_1)$. Thus
\begin{equation}
         \int_{s(m)}^{s(m+1)}  \big(\Delta x_m\big)^{\ell} dx
        =  \left[ \frac{(\Delta x_m)^{\ell+1}}{\ell+1} 
           \right]^{s(m+1)}_{s(m)} = \frac{V^{-\ell-1}
            }{\ell+1}.
\end{equation}
\vskip 0.3cm

It then can be verified that
\begin{eqnarray}
F(V)&=&\sum_{m=m_1}^{m_2}\sum_{n=0}^{m} f(m,n)e^{Vg(m,n)}
\nonumber\\
&\times&
\sum_{\ell=0}^\infty \frac{V^{\ell+2}}{\ell!} \int_{s(m)}^{s(m+1)}\int_{s(n)}^{s(n+1)}\Big[\big(g_x(m,n)\Delta x_m+g_y(m,n)\Delta y_m\big)^{\ell}\Big] dydx
\nonumber\\[5pt]
  &=&\sum_{m=m_1}^{m_2}\sum_{n=0}^{m} f(m,n)e^{Vg(m,n)}
\nonumber\\
&\times&
\sum_{\ell=0}^\infty \frac{V^{\ell+2}}{\ell!} \int_{0}^{\frac{1}{V}}\int_{0}^{\frac{1}{V}}\Big[\big(g_x(m,n)x
+g_y(m,n) y\big)^{\ell}\Big] dydx
\nonumber\\
&=&\sum_{m=m_1}^{m_2}\sum_{n=0}^{m} \frac{f(m,n)e^{Vg(m,n)}}{\ell!} \sum_{\ell=0}^\infty  V^{\ell+2}\int_{0}^{\frac{1}{V}}\int_{0}^{\frac{1}{V}}\sum_{k=0}^{\ell}
{\ell\choose k}
\nonumber\\[5pt]
&& \hspace{0.4in} \times\Big[g_x^k(m,n)x^k
g_y^{\ell-k}(m,n)y^{\ell-k}\Big]dydx
\nonumber\\
&=&\sum_{m=m_1}^{m_2}\sum_{n=0}^{m} f(m,n)e^{Vg(m,n)} \sum_{\ell=0}^\infty\sum_{k=0}^{\ell}\frac{g_x^k(m,n)g_y^{\ell-k}(m,n)}{(k+1)!(\ell-k+1)!}
\nonumber\\[5pt]
&=&\sum_{m=m_1}^{m_2}\sum_{n=0}^{m} \frac{f(m,n)e^{Vg(m,n)}}{g_x(m,n)g_y(m,n)}\sum_{\ell=0}^\infty\sum_{k=1}^{\ell+1}\frac{g_x^k(m,n)g_y^{\ell-k+2}(m,n)}{k!(\ell-k+2)!}
\nonumber
\end{eqnarray}

\begin{eqnarray}
&=&\sum_{m=m_1}^{m_2}\sum_{n=0}^{m}  \frac{f(m,n)e^{Vg(m,n)}}{g_x(m,n)g_y(m,n)}\sum_{\ell=0}^\infty\frac{1}{(\ell+2)!}
\nonumber\\[5pt]
&& \hspace{0.4in}
 \times\Big[(g_x(m,n)+g_y(m,n))^{\ell+2}-g_x^{\ell+2}(m,n)-g_y^{\ell+2}(m,n)\Big]
\nonumber\\[5pt]
&=&\sum_{m=m_1}^{m_2}\sum_{n=0}^{m}  
  \frac{f(m,n)e^{Vg(m,n)}}{g_x(m,n)g_y(m,n)}
  \Big[e^{g_x(m,n)+g_y(m,n)}-e^{g_x(m,n)}
\nonumber\\[5pt]
&& \hspace{0.4in} -e^{g_y(m,n)}+1\Big]\cr
&=&\sum_{m=m_1}^{m_2}\sum_{n=0}^{m} f(m,n)e^{Vg(m,n)} \left[ \frac{(e^{g_x(m,n)}-1)(e^{g_y(m,n)}-1)}
{g_x(m,n)g_y(m,n)}\right].
\label{eq00065}
\end{eqnarray}
Combining Eqs. (\ref{def_FV}) and (\ref{eq00065}) we obtain
\begin{equation}
     \textrm{Eq. } (\ref{eq00065}) = F(V) 
           = \left(\int_a^b\int_0^x  f(x,y)
                e^{Vg(x,y)}dydx\right) 
                \Big(V^2+ O(V)\Big). 
\end{equation}

Therefore, letting 
\[
   \tilde{f}(m,n) = f(m,n)\left[ \frac{(e^{g_x(m,n)}-1)(e^{g_y(m,n)}-1)}
{g_x(m,n)g_y(m,n)}\right],
\]
we have
\begin{equation}
 \sum_{m=m_1}^{m_2}\sum_{n=0}^{m} \tilde{f}(m,n)e^{Vg(m,n)}
   = \left(\int_a^b\int_0^x\frac{\tilde{f}(x,y)\ g_x\ 
             g_y\ e^{Vg(x,y)}}
    {(e^{g_x}-1)(e^{g_y}-1)}dydx\right)\Big(V^2+ O(V)\Big).
\label{the_result00}
\end{equation}

To show Eq. (\ref{T12_exact}) to Eq. (\ref{asym_4_T}), we use
$a=x$, $b=x_2$,
\[
g(x,y)=\phi_0(x)-\phi_0(y)=\int_y^x\ln
        \left(\frac{\lambda_0(z)}{\mu_0(z)}\right) dz,
         \textrm{ and } \ 
\]
\[
\tilde{f}(x,y)e^{Vg(x,y)}=\exp\Big[
      V\big(\Phi(x,V)-\Phi(y,V)\big)\Big].
\]

\subsection{Laplace's method for integrals beyond leading order}

One can find the materials in this section from
classic texts such as \cite{murray,omalley}.  They are 
included here for the convenience of the reader.  
We first have two important formulae:
\begin{equation}
	\frac{2}{\sqrt{\pi}} \int_0^x e^{-t^2} dt = \textrm{erf}(x),
\label{i_eq_1}
\end{equation}
and for large $x$
\begin{equation}
	\textrm{erf}(x) = 1-\frac{e^{-x^2}}{x\sqrt{\pi}}
		\left[1+\sum_{n=1}^{\infty} (-1)^n
			\frac{1\cdot 3\cdot 5\ ...\ (2n-1)}{(2x^2)^n}
		\right].
\label{i_eq_2}
\end{equation}
Using Eqs. (\ref{i_eq_1}) and (\ref{i_eq_2}), we can evaluate
the following integral for $z > 0$:
\begin{eqnarray}
	G(x) &=& \int_0^x e^{-\frac{(y-z)^2}{\epsilon}}dy
\label{the_LI}\\ 
	&=& \frac{\sqrt{\pi}}{2}\left[
		\textrm{erf}\left(\frac{x-z}{\sqrt{\epsilon}}\right)+
		 \textrm{erf}\left(\frac{z}{\sqrt{\epsilon}}\right) \right]
\nonumber\\[7pt]
	&\approx& \left\{\begin{array}{ccc}
		\frac{\epsilon}{2} \left[
	\frac{e^{-(x-z)^2/\epsilon}}{z-x}
		\left(1-\frac{\epsilon}{2(x-z)^2}\right)
	-\frac{e^{-z^2/\epsilon}}{z}\left(1-\frac{\epsilon}{2z^2} \right)	
	 \right] & & x < z 
\\[9pt]
		\frac{\sqrt{\pi\epsilon}}{2}-
	\frac{\epsilon e^{-z^2/\epsilon}}{2z}\left(1-\frac{\epsilon}{2z^2} \right)
			 & & x = z 
\\[9pt]
			\sqrt{\pi\epsilon}-\frac{\epsilon}{2} \left[
	\frac{e^{-z^2/\epsilon}}{z}\left(1-\frac{\epsilon}{2z^2} \right)	
	+ \frac{e^{-(x-z)^2/\epsilon}}{x-z}
		\left(1-\frac{\epsilon}{2(x-z)^2}\right) \right] & & x > z  \end{array}
		\right.
\nonumber\\
\end{eqnarray}
Near $x=z$, there is a boundary layer.  Let $x\mapsto\sigma=(x-z)/\sqrt{\epsilon}$ 
and $z\mapsto\eta=(y-z)/\sqrt{\epsilon}$, then the integral in (\ref{the_LI}) has
an inner expression
\begin{eqnarray}
	\widetilde{G}(\sigma) = 
	\int_0^x e^{-\frac{(y-z)^2}{\epsilon}}dy  &=&
	\sqrt{\epsilon}\int_{-\sigma}^{z/\sqrt{\epsilon}} e^{-\eta^2}d\eta
\nonumber\\
	&=& \frac{\sqrt{\epsilon\pi}}{2}\left[
		\textrm{erf}\left(\frac{z}{\sqrt{\epsilon}}\right)
		+\textrm{erf}\left(\sigma\right)\right]
\nonumber\\
	&\approx& \frac{\sqrt{\epsilon\pi}}{2}\left(1
		+\textrm{erf}\left(\sigma\right)\right)
		-\frac{\epsilon e^{-z^2/\epsilon}}{2z}
		\left[1-\frac{\epsilon}{(2z^2)}
		\right].
\nonumber\\[-9pt]
\end{eqnarray}	
We have a matched asymptotics: $\widetilde{G}(0) = G(z)$,
\begin{eqnarray}
	\lim_{\sigma\rightarrow -\infty} \widetilde{G}(\sigma)
		&=& \frac{\sqrt{\epsilon} e^{-\sigma^2}}
		{2|\sigma|} \left(1-\frac{1}{2\sigma^2}\right)  
		= \lim_{x\rightarrow z^-} G(x), 
\\ [9pt]
	\lim_{\sigma\rightarrow \infty} \widetilde{G}(\sigma)
		&=& \sqrt{\sigma\pi}-\frac{\sqrt{\epsilon} e^{-\sigma^2}}
		{2\sigma} \left(1-\frac{1}{2\sigma^2}\right)   
		=\lim_{x\rightarrow z^+} G(x).
\end{eqnarray}

Therefore, for a $\phi(y)$ with a unique minimum at 
$y=x^{\ddag}>0$:
\begin{eqnarray}
	&& \frac{1}{V}\ln\int_0^x e^{-V\phi(y)}dy 
\nonumber\\[5pt]
     &=&  \left\{\begin{array}{lcc}
		-\phi(x)+\frac{1}{V}\ln
	   	\int_0^x e^{-V\left[\phi'(x)(y-x)
		+\frac{1}{2}\phi''(x)(y-x)^2\right]}
		dy & & x < x^{\ddag}		
\\[9pt]
		-\phi(x^{\ddag})+\frac{1}{V}\ln
	   \int_0^x e^{-V\left[\frac{1}{2}\phi''(x^{\ddag})(y-x^{\ddag})^2\right]}
		dy    && x > x^{\ddag}
		\end{array}\right.
\nonumber\\[9pt]
	&=&  \left\{\begin{array}{lcc}
	 -\phi(x)-\frac{1}{V}\ln(V|\phi'(x)|) 
	 - \frac{1}{V^2}\frac{\phi''(x)}{\left(\phi'(x)\right)^2} & & x<x^{\ddag} 
\\[9pt]
	-\phi(x^{\ddag})+\frac{1}{2V}\ln\frac{2\pi}{V\phi''(x^{\ddag})}	
		-\frac{1}{V\sqrt{2\pi V\phi''(x^{\ddag})}} \left[
	 \frac{e^{-V\phi''(x^{\ddag})\left(x-x^{\ddag}\right)^2/2}}{x-x^{\ddag}} \right] & & x > x^{\ddag}  \end{array}
		\right.
\nonumber
\end{eqnarray}
And within the boundary layer $x\in\left(x^{\ddag}-\sqrt{\epsilon},x^{\ddag}+\sqrt{\epsilon}\right)$
where $\epsilon=\sqrt{2/\big(V\phi''(x^{\ddag})\big)}$, we have an inner expansion:
\begin{equation}
 \frac{1}{V}\ln\left[1
	+\textrm{erf}\left(\frac{\sqrt{V\phi''(x^{\ddag})}}{2}(x-x^{\ddag})\right)\right].
\label{sto_pot}
\end{equation}
With increasing $x$, the stochastic potential in 
Eq. (\ref{sto_pot}) switches from an ``enthalpic dominant''
\[
      -\phi(x)-\frac{1}{V}\ln(V|\phi'(x))
	 + o\left(V^{-1}\right), \ \ x< x^{\ddag},
\] 
to a constant independent of $x$
\[   
	 	-\phi(x^{\ddag})
		+\frac{1}{2V}\ln\frac{2\pi}{V\phi''(x^{\ddag})}, 
		\ \ \
		x > x^{\ddag}.
\]
The enthalpic term is in fact the limit
\begin{equation}
	\lim_{V\rightarrow\infty} \frac{1}{V}
	\ln\int_0^x e^{-V\phi(y)}dy 
		= -\inf_{y\in[0,x]} \phi(y).
\end{equation}

	Now for the cases in which $\phi(y)$ has a single maximum at
$y=x^{\ddag}>0$, we have:
\[
   \frac{1}{V}\ln\int_0^x e^{-V\phi(y)}dy = \hspace{3.5in}
\]
\begin{equation}
	 \left\{\begin{array}{lcc}
		-\phi(0)-\frac{1}{V}\ln\big(V\phi'(0)\big)-\frac{1}{V}
	 \frac{\phi''(0)}{\left(\phi'(0)\right)^2} && 0<x<x^{\ddag} 
\\[9pt]

		-\phi(0)+\frac{1}{V}\ln\left[\frac{1}{V\phi'(0)}
			+ \frac{1}{V|\phi'(x)|}e^{-V\big(\phi(x)-\phi(0)\big)}
			\right]& & x > x^{\ddag}  \end{array}
		\right.
\label{maximum}
\end{equation}
In this case, there is a switching from ``entropic dominance'' 
constant when $x<x^{\ddag}<x^*$ where $\phi(x^*)=\phi(0)$ to 
``enthalpic dominance'' $-\phi(x)$ when $x>x^*$.

\subsection{Detailed derivation for Eq. (\ref{new_k_eq})}
\label{sec.vii_c}

For the sake of convenience, we rewrite Eq. (\ref{asym_4_T}) as follows
\begin{equation}
V\int_a^b\int_0^x\hat f(x,y)e^{Vg(x,y)}dydx
\end{equation}
where, the same as in Sec. (\ref{sec.vii_b}), $g(x,y)=\phi_0(x)-\phi_0(y)=\int_y^x\ln\frac{\lambda_0(z)}{\mu_0(z)}dz$.
If there exists a point $(x^*, y^*)$ in $\{(x,y)|x_1^*<x<x_2^*, 0<y<x\}$ that satisfies $g(x,y)\le g(x^*, y^*)$, i.e.
\begin{eqnarray}
\label{assumption}
&&g_x(x^*, y^*)=g_y(x^*, y^*)=0, 
\ \ \ g_{xx}(x^*, y^*)<0, g_{yy}(x^*, y^*)<0, 
\nonumber\\[5pt]
  &&g_{xy}^2(x^*, y^*)-g_{xx}(x^*, y^*)g_{yy}(x^*, y^*)< 0,
\end{eqnarray}
and
\begin{eqnarray}
\hat f(x^*, y^*)\ne 0,
\end{eqnarray}
then for large $V$, denoting $\Delta x=(x-x^*)$ and 
$\Delta y=(y-y^*)$, 
\begin{eqnarray}
&&V\int_a^b\int_0^x\hat f(x,y)e^{Vg(x,y)}dydx
\nonumber\\
&\approx&
V\hat f(x^*, y^*)e^{Vg(x^*, y^*)} \int_a^b\int_0^xe^{V\left[\frac12 g_{xx}(\Delta x)^2+\frac12 g_{yy}(\Delta y)^2+g_{xy}(\Delta x)(\Delta y)\right]}dydx
\nonumber\\
&\approx&
V\hat f(x^*, y^*)e^{Vg(x^*, y^*)} \int_{x^*-\epsilon}^{x^*+\epsilon}\int_{y^*-\epsilon}^{y^*+\epsilon}
e^{V\left[\frac12 g_{xx}(\Delta x)^2+\frac12 g_{yy}(\Delta y)^2+g_{xy}(\Delta x)(\Delta y)\right]}dydx
\nonumber
\end{eqnarray}
\begin{eqnarray}
&\approx&
V\hat f(x^*, y^*)e^{Vg(x^*, y^*)} \int_{-\infty}^{\infty}\int_{-\infty}^{\infty}
e^{V\left[\frac12 g_{xx}(\Delta x)^2+\frac12 g_{yy}(\Delta y)^2+g_{xy}(\Delta x)(\Delta y)\right]}dydx
\nonumber\\
&\approx&
V\hat f(x^*, y^*)e^{Vg(x^*, y^*)} \int_{-\infty}^{\infty}\int_{-\infty}^{\infty}
e^{V\left[\frac12 g_{xx}x^2+\frac12 g_{yy}x^2+g_{xy}xy\right]}dydx
\nonumber\\
&=&\frac{2\pi \hat f(x^*, y^*)e^{Vg(x^*, y^*)}}{\sqrt{g_{xx}(x^*, y^*)g_{yy}(x^*, y^*)-g_{xy}^2(x^*, y^*)}}.
\label{derivation 34}
\end{eqnarray}
So, from Eqs. (\ref{asym_4_T}) and (\ref{derivation 34}), 
one can find
\begin{eqnarray}
  && T_{x_1^*\rightarrow x^*_2} 
\nonumber\\[5pt]
 &\approx&
\frac{\ln\frac{\lambda_0(x^*)}{\mu_0(x^*)}\ln\frac{\mu_0(y^*)}{\lambda_0(y^*)}}{\left(\frac{\lambda_0(x^*)}{\mu_0(x^*)}-1\right)\left(\frac{\mu_0(y^*)}{\lambda_0(y^*)}-1\right)}
  \frac{2\pi e^{V[\phi_0(x^*,V)-\phi_0(y^*,V)]+\phi_1(x^*,V)-\phi_1(y^*,V)}}{\lambda_0(x^*)\sqrt{g_{xx}(x^*, y^*)g_{yy}(x^*, y^*)-g_{xy}^2(x^*, y^*)}}\cr
  &=&
\frac{\ln\frac{\lambda_0(x^*)}{\mu_0(x^*)}\ln\frac{\mu_0(y^*)}{\lambda_0(y^*)}}{\left(\frac{\lambda_0(x^*)}{\mu_0(x^*)}-1\right)\left(\frac{\mu_0(y^*)}{\lambda_0(y^*)}-1\right)}
  \frac{2\pi e^{V[\phi_0(x^*,V)-\phi_0(y^*,V)]+\phi_1(x^*,V)-\phi_1(y^*,V)}}{\lambda_0(x^*)
  \sqrt{\left(\frac{\lambda_0'(x^*)}{\lambda_0(x^*)}-\frac{\mu_0'(x^*)}{\mu_0(x^*)}\right)
  \left(\frac{\mu_0'(y^*)}{\mu_0(y^*)}-\frac{\lambda_0'(y^*)}{\lambda_0(y^*)}\right)}}\cr
  &=&
  \frac{2\pi e^{V[\phi_0(x^*,V)-\phi_0(y^*,V)]+\phi_1(x^*,V)-\phi_1(y^*,V)}}{\lambda_0(x^*)
  \sqrt{\left(\frac{\lambda_0'(x^*)}{\lambda_0(x^*)}-\frac{\mu_0'(x^*)}{\mu_0(x^*)}\right)
  \left(\frac{\mu_0'(y^*)}{\mu_0(y^*)}-\frac{\lambda_0'(y^*)}{\lambda_0(y^*)}\right)}}\cr
    &=&
  \frac{2\pi e^{V[\Phi_0(x^*,V)-\Phi_0(y^*,V)]}}{\lambda_0(x^*)
  \sqrt{-\phi_0''(x^*)\phi_0''(y^*)}},
\end{eqnarray}
where the third equation is because that, at the maximum point $(x^*, y^*)$, $\lambda_0(x^*)=\mu_0(x^*)$, $\lambda_0(y^*)=\mu_0(y^*)$ and
$\lim_{z\to 1}\frac{\ln z}{z-1}=1$. The last equation is because
Eq. (\ref{assumption}) implies $\phi_0''(x^*)\phi_0''(y^*)<0$.

\end{document}